\documentclass[12pt]{amsart}

\usepackage{times,epsf, rlepsf,amsmath, hyperref}

\begin{document}

\newtheorem{thm}{Theorem}[section]
\newtheorem{lem}[thm]{Lemma}
\newtheorem{cor}[thm]{Corollary}

\theoremstyle{definition}
\newtheorem{defn}{Definition}[section]

\theoremstyle{remark}
\newtheorem{rmk}{Remark}[section]

\def\square{\hfill${\vcenter{\vbox{\hrule height.4pt \hbox{\vrule
width.4pt height7pt \kern7pt \vrule width.4pt} \hrule height.4pt}}}$}

\def\T{\mathcal T}

\newenvironment{pf}{{\it Proof:}\quad}{\square \vskip 12pt}

\title{On the number of solutions to the Asymptotic Plateau Problem}
\author{Baris Coskunuzer}
\address{Department of Mathematics \\ Yale University \\ New Haven, CT 06520}
\curraddr{Koc University, Department of Mathematics, Istanbul, Turkey}

\email{bcoskunuzer@ku.edu.tr}
\thanks{The author is supported by NSF Grant DMS-0603532, EU-FP7 Grant IRG-226062 and TUBITAK Grant 107T642}

\maketitle


\newcommand{\cirD}{\overset{\circ}{D}}
\newcommand{\Si}{S^2_{\infty}({\mathbf H}^3)}
\newcommand{\SI}{S^{n-1}_{\infty}({\mathbf H}^n)}
\newcommand{\PI}{\partial_{\infty}}
\newcommand{\BH}{\mathbf H}
\newcommand{\BR}{\mathbf R}
\newcommand{\BC}{\mathbf C}
\newcommand{\BZ}{\mathbf Z}
\newcommand{\BHH}{{\mathbf H}^3}

\begin{abstract}

We give a simple topological argument to show that the number of solutions to the asymptotic Plateau problem in
hyperbolic space is generically unique. In particular, we show that the space of closed, codimension-$1$ submanifolds
of $\SI$ which bound a unique absolutely area minimizing hypersurface in $\BH^n$ is dense in the space of
codimension-$1$ closed submanifolds of $\SI$. In dimension $3$, we also prove that the set of uniqueness curves in
$\Si$ for least area planes is generic in the set of Jordan curves in $\Si$. We also give some nonuniqueness results
for dimension $3$, too.

\end{abstract}

\section{Introduction}

The asymptotic Plateau problem in hyperbolic space asks the existence of an absolutely area minimizing hypersurface
$\Sigma\subset\BH^n$ asymptotic to given closed codimension-$1$ submanifold $\Gamma$ in $\SI$. This problem is solved
by Michael Anderson in his seminal paper \cite{A1}. He proved the existence of a solution for any given closed
submanifold in the sphere at infinity. Then, Hardt and Lin studied the asymptotic regularity of these solutions in
\cite{HL}, \cite{Li}. Lang generalized Anderson's methods to solve this problem in Gromov-Hadamard spaces in \cite{L}.

On the other hand, on the number of the absolutely area minimizing hypersurfaces for a given asymptotic boundary, there
are a few results so far. In \cite{A1}, Anderson showed that if the given asymptotic boundary $\Gamma$ bounds a convex
domain in $\SI$, then there exists a unique absolutely area minimizing hypersurface in $\BH^n$ . Then, Hardt and Lin
generalized this result to the closed codimension-$1$ submanifolds bounding star shaped domains in $\SI$ in \cite{HL}.
Recently, the author showed a generic uniqueness result in dimension $3$ for least area planes in \cite{Co1},
\cite{Co2}. For other results on asymptotic Plateau problem, see the survey article \cite{Co5}.

In this paper, we will show that the space of closed codimension-$1$ submanifolds in $\SI$ bounding a unique absolutely
area minimizing hypersurface in $\BH^n$ is dense in the space of all closed codimension-$1$ submanifolds in $\SI$ by
using a simple topological argument. On the other hand, in dimension $3$, we show that the space of Jordan curves in
$\Si$ bounding a unique area minimizing surface in $\BH^3$ is generic in the space of Jordan curves in $\Si$.
Similarly, we show that the same is true for the least area planes in $\BH^3$, too.

The main results of the paper are the following. The first one is about least area planes in $\BH^3$.

\vspace{0.3cm}

\noindent \textbf{Theorem 3.3.} Let $A$ be the space of simple closed curves in $\Si$ and let $A'\subset A$ be the
subspace containing the simple closed curves in $\Si$ bounding a unique least area plane in $\BH^3$. Then, $A'$ is
generic in $A$, i.e. $A-A'$ is a set of first category.

\vspace{0.3cm}

The second result is about absolutely area minimizing hypersurfaces.

\vspace{0.3cm}

\noindent \textbf{Theorem 4.4.} Let $B$ be the space of connected closed codimension-$1$ submanifolds of $\SI$, and let
$B'\subset B$ be the subspace containing the closed submanifolds of $\SI$ bounding a unique absolutely area minimizing
hypersurface in $\BH^n$. Then $B'$ is dense in $B$.

\vspace{0.3cm}

Indeed, this subspace of closed submanifolds in $\SI$ bounding a unique absolutely area minimizing hypersurface is not
only dense, but also generic in some sense. See Remark 4.2. On the other hand, in dimension $3$, by combination of
these two theorems, we get the following corollary.

\vspace{0.3cm}

\noindent \textbf{Corollary 4.5.} Let $A$ be the space of simple closed curves in $\Si$ and let $A'\subset A$ be the
subspace containing the simple closed curves in $\Si$ bounding a unique absolutely area minimizing surface in $\BH^3$.
Then, $A'$ is generic in $A$, i.e. $A-A'$ is a set of first category.

\vspace{0.3cm}

The short outline of the technique for generic uniqueness is the following: For simplicity, we will focus on the case
of the least area planes in $\BH^3$. Let $\Gamma_0$ be a simple closed curve in $\Si$. First, we will show that either
there exists a unique least area plane $\Sigma_0$ in $\BH^3$ with $\PI\Sigma_0=\Gamma_0$, or there exist two {\em
disjoint} least area planes $\Sigma_0^+ , \Sigma_0^-$ in $\BH^3$ with $\PI\Sigma_0^\pm=\Gamma_0$. Now, take a small
neighborhood $N(\Gamma_0)\subset \Si$ which is an annulus. Then foliate $N(\Gamma_0)$ by simple closed curves
$\{\Gamma_t\}$ where $t\in(-\epsilon, \epsilon)$, i.e. $N(\Gamma_0) \simeq \Gamma\times (-\epsilon, \epsilon)$. By the
above fact, for any $\Gamma_t$ either there exists a unique least area plane $\Sigma_t$, or there are two least area
planes $\Sigma_t^\pm$ disjoint from each other. Also, since these are least area planes, if they have disjoint
asymptotic boundary, then they are disjoint by Meeks-Yau exchange roundoff trick. This means, if $t_1<t_2$, then
$\Sigma_{t_1}$ is disjoint and \textit{below} from $\Sigma_{t_2}$ in $\BH^3$. Consider this collection of least area
planes. Note that for curves $\Gamma_t$ bounding more than one least area plane, we have a canonical region $N_t$ in
$\BH^3$ between the disjoint least area planes $\Sigma_t^\pm$, see Figure 1.

Now, $N(\Gamma)$ separates $\Si$ into two parts, and take a geodesic $\beta\subset \BH^3$ which is asymptotic to two
points belongs to these two different parts. This geodesic is transverse to the collection of these least area planes
asymptotic to the curves in $\{\Gamma_t\}$. Also, a finite segment of this geodesic intersects the entire collection.
Let the length of this finite segment be $C$. Now, the idea is to consider the {\em thickness} of the neighborhoods
$N_t$ assigned to the asymptotic curves $\{\Gamma_t\}$. Let $s_t$ be the length of the segment $I_t$ of $\beta$ between
$\Sigma_t^+$ and $\Sigma_t^-$, which is the {\em width} of $N_t$ assigned to $\Gamma_t$. Then, the curves $\Gamma_t$
bounding more than one least area plane have positive width, and contributes to total thickness of the collection, and
the curves bounding unique least area plane has $0$ width and do not contribute to the total thickness. Since
$\sum_{t\in(-\epsilon, \epsilon)} s_t < C$, the total thickness is finite. This implies for only countably many
$t\in(-\epsilon, \epsilon)$, $s_t>0$, i.e. $\Gamma_t$ bounds more than one least area plane. For the remaining
uncountably many $t\in(-\epsilon, \epsilon)$, $s_t=0$, and there exists a unique least area plane for those $t$. This
proves the space of Jordan curves of uniqueness is dense in the space of Jordan curves in $\Si$. Then, we will show
this space is not only dense, but also generic.

We should note that this technique is quite general and it can be applied to many different settings of the Plateau
problem (See Concluding Remarks).

On the other hand, after the above generic uniqueness results, it is a reasonable question whether all simple closed
curves in $\Si$ bound a unique absolutely area minimizing surface or a unique least area plane. It is still not known
whether all simple closed curves in $\Si$ have a unique solution to the asymptotic Plateau problem or not. The only
known results about nonuniqueness also come from Anderson in \cite{A2}. He constructs examples of simple closed curves
in $\Si$ bounding more than one complete \textit{minimal surface} in $\BH^3$. These examples are also area minimizing
in their topological class. However, none of them are absolutely area minimizing, i.e. a solution to the asymptotic
Plateau problem.

In this paper, we show the existence of simple closed curves in $\Si$ with nonunique solution to the asymptotic Plateau
problem. In other words, we show that there are examples of simple closed curves in $\Si$ which are the asymptotic
boundaries of more than one absolutely area minimizing surface.

\vspace{0.3cm}

\noindent {\bf Theorem 5.2.} There exists a simple closed curve $\Gamma$ in $\Si$ such that $\Gamma$ bounds more than
one absolutely area minimizing surface $\{\Sigma_i\}$ in $\BHH$, i.e. $\PI \Sigma_i = \Gamma$.

\vspace{0.3cm}

Also, by using similar ideas, we show nonuniqueness for least area planes case in dimension $3$.

\vspace{0.3cm}

\noindent {\bf Theorem 5.3.} There exists a simple closed curve $\Gamma$ in $\Si$ such that $\Gamma$ bounds more than
one least area plane $\{P_i\}$ in $\BHH$, i.e. $\PI P_i = \Gamma$.

\vspace{0.3cm}

The organization of the paper is as follows: In the next section we will cover some basic results which will be used in
the following sections. In section 3, we will show the genericity result for least area planes in $\BH^3$. Then in
section 4, we will show the results on absolutely area minimizing hypersurfaces in $\BH^n$. In Section 5, we will prove
the nonuniqueness results. Finally in section 6, we will have some concluding remarks.

\subsection{Acknowledgements:}

I would like to thank Yair Minsky and Michael Anderson for very useful conversations.

\section{Preliminaries}

In this section, we will overview the basic results which we will use in the following sections. For details on the
notions and results in this section, see the survey article \cite{Co5}.

First, we will give the definitions of area minimizing surfaces. First set of the definitions are about compact
surfaces and hypersurfaces. The second set of the definitions are their generalizations to the noncompact surfaces and
hypersurfaces.

\begin{defn}(Compact Case) A {\em least area disk} (area minimizing disk) is a disk which has the smallest area among the disks with the same boundary.
An {\em absolutely area minimizing surface} is a surface which has the smallest area among all the surfaces (with no
topological restriction) with the same boundary. An {\em absolutely area minimizing hypersurface} is a hypersurface
which has the smallest volume among all hypersurfaces with the same boundary.
\end{defn}

\begin{defn}(Noncompact Case) A {\em least area plane} is a plane such that any compact subdisk in the plane
is a least area disk. We will also call a complete noncompact surface as {\em absolutely area minimizing surface} if
any compact subsurface is an absolutely area minimizing surface. Similarly, we will call a complete noncompact
hypersurface as {\em absolutely area minimizing hypersurface}, if any compact part (codimension-$0$ submanifold with
boundary) of the hypersurface is an absolutely area minimizing hypersurface.
\end{defn}

Now, we will quote the basic results on asymptotic Plateau problem.

\begin{lem} \cite{A1} Let $\Gamma$ be a codimension-$1$ closed submanifold of $\SI$. Then there exists a complete, absolutely
area minimizing $n-1$-rectifiable current $\Sigma$ in $\BH^n$ with $\PI \Sigma=\Gamma$.
\end{lem}

Note that the rectifiable current here is indeed a smooth hypersurface of $\BH^n$ except for a singular set of
Hausdorff dimension at most $n-8$ by the regularity result stated below. For convenience, we will call area minimizing
codimension-$1$ rectifiable currents as area minimizing hypersurfaces throughout the paper.

\begin{lem} \cite{A2} Let $\Gamma$ be a Jordan curve in $\Si$. Then, there exists a complete, least area plane
$\Sigma$ in $\BH^3$ asymptotic to $\Gamma$.
\end{lem}

The following fact about interior regularity theory of geometric measure theory is well-known.

\begin{lem} \cite{Fe} Let $\Sigma$ be a $(n-1)$-dimensional area minimizing rectifiable current. Then $\Sigma$ is a smooth,
embedded manifold in the interior except for a singular set of Hausdorff dimension at most $n-8$.
\end{lem}

Finally, we will state a theorem about limits of sequences of least area planes. Here, the limit is the pointwise limit
of the planes, and for each limit point, there is a disk containing the point in the limit set such that it is the
limit of a sequence of subdisks in the planes [Ga, Lemma 3.3].

\begin{lem}\cite{Ga}
Let $\{\Sigma_i\}$ be a sequence of least area planes in $\BH^3$ with $\PI \Sigma_i = \Gamma_i \subset \Si$ simple
closed curve for any $i$. If $\Gamma_i\rightarrow\Gamma$, then there exists a subsequence $\{\Sigma_{i_j}\}$ of
$\{\Sigma_i\}$ such that $\Sigma_{i_j}\rightarrow\widehat{\Sigma}$ a collection of least area planes whose asymptotic
boundaries are $\Gamma$.
\end{lem}

The next theorem is a similar limit theorem about absolutely area minimizing hypersurfaces in $\BH^n$.

\begin{lem} Let $\{\Gamma_i\}$ be a sequence of connected closed codimension-$1$ submanifolds in $\SI$ which are pairwise disjoint.
Let $\{\Sigma_i\}$ be a sequence of complete absolutely area minimizing hypersurfaces in $\BH^n$ with $\PI(\Sigma_i) =
\Gamma_i$. If $\Gamma_i$ converges to a closed codimension-$1$ submanifold $\Gamma$ in $\SI$, then there exists a
subsequence of $\{\Sigma_i\}$ which converges to a complete absolutely area minimizing hypersurfaces $\Sigma$ in
$\BH^n$ with $\PI \Sigma = \Gamma$.
\end{lem}

\begin{pf} By the proof of the Lemma 2.1, all we need to show that we can induce a suitable sequence of
compact absolutely area minimizing hypersurfaces $\{S_i\}$ from $\{\Sigma_i\}$ where $\gamma_i=\partial S_i$
converges to $\Gamma$ in $\Si$.

Let $K>0$ be sufficiently large so that $\Sigma_i \cap N_K(CH(\Gamma)) \neq \emptyset$ where $N_K(CH(\Gamma))$
is the $K$ neighborhood of the convex hull of $\Gamma$. Then, let $S_i = N_K(CH(\Gamma)) \cap \Sigma_i$.

Since $\Sigma_i \subset CH(\Gamma_i)$ and $\Gamma_i \rightarrow \Gamma$, then $\partial S_i$ must converge to $\Gamma$
asymptotically. Then, by the proof of Lemma 2.1 with slight modification for the uniform mass bounds coming from
$N_K(CH(\Gamma))$, there exists a subsequence $\{S_{i_j}\}$ which converges to a complete absolutely area minimizing
hypersurface $\Sigma$ in $\BH^n$ with $\PI \Sigma = \Gamma$.
\end{pf}

\section{Least Area Planes in $\BH^3$}

In this section, we will prove that the space of simple closed curves in $\Si$ bounding a unique least area plane in
$\BH^3$ is generic in the space of simple closed curves in $\Si$.

First, we will show that if two least area planes have disjoint asymptotic boundaries, then they are disjoint.

\begin{lem}
Let $\Gamma_1$ and $\Gamma_2$ be two disjoint simple closed curves in $\Si$. If $\Sigma_1$ and $\Sigma_2$ are least
area planes in $\BH^3$ with $\PI \Sigma_i = \Gamma_i$, then $\Sigma_1$ and $\Sigma_2$ are disjoint, too.
\end{lem}

\begin{pf}
Assume that $\Sigma_1\cap\Sigma_2\neq\emptyset$. Since asymptotic boundaries $\Gamma_1$ and $\Gamma_2$ are
disjoint, the intersection cannot contain an infinite line. So, the intersection between $\Sigma_1$ and
$\Sigma_2$ must contain a simple closed curve $\gamma$. Since $\Sigma_1$ and $\Sigma_2$ are also minimal, the
intersection must be transverse on a subarc of $\gamma$ by maximum principle.

Now, $\gamma$ bounds two area minimizing disks $D_1$ and $D_2$ in $\BH^3$, with $D_i\subset\Sigma_i$. Now, take a
larger subdisk $E_1$ of $\Sigma_1$ containing $D_1$, i.e. $D_1\subset E_1 \subset \Sigma_1$. By definition, $E_1$ is
also an area minimizing disk. Now, modify $E_1$ by swaping the disks $D_1$ and $D_2$. Then, we get a new disk $E_1 '=
\{E_1 - D_1\} \cup D_2$. Now, $E_1$ and $E_1 '$ have same area, but $E_1 '$ have a folding curve along $\gamma$. By
smoothing out this curve as in \cite{MY}, we get a disk with smaller area, which contradicts to $E_1$ being area
minimizing. Note that this technique is known as Meeks-Yau exchange roundoff trick.
\end{pf}

The following lemma is very essential for our technique. Mainly, the lemma says that for any given simple closed curve
$\Gamma$ in $\Si$, either there exists a unique least area plane $\Sigma$ in $\BH^3$ asymptotic to $\Gamma$, or there
exist two least area planes $\Sigma^\pm$ in $\BH^3$ which are asymptotic to $\Gamma$ and disjoint from each other. Even
though this lemma is also proven in \cite{Co4}, because of its importance for the technique, and to set the notation
for the main result, we give a proof here. Note that Brian White proved a similar version of this lemma by using
geometric measure theory methods in \cite{Wh1}.

\begin{lem}
Let $\Gamma$ be a simple closed curve in $\Si$. Then either there exists a unique least area plane $\Sigma$ in $\BH^3$
with $\PI\Sigma=\Gamma$, or there are two canonical disjoint extremal least area planes $\Sigma^+$ and $\Sigma^-$ in
$\BH^3$ with $\PI \Sigma^\pm = \Gamma$. Moreover, any least area plane $\Sigma '$ with $\PI\Sigma '= \Gamma$ is
disjoint from $\Sigma^\pm$, and it is captured in the region bounded by $\Sigma^+$ and $\Sigma^-$ in $\BH^3$.
\end{lem}

\begin{pf}
Let $\Gamma$ be a simple closed curve in $\Si$. $\Gamma$ separates $\Si$ into two parts, say $D^+$ and $D^-$.
Define sequences of pairwise disjoint simple closed curves $\{\Gamma_i^+\}$ and $\{\Gamma_i^-\}$ such that
$\Gamma_i^+\subset D^+$, and $\Gamma_i^-\subset D^-$ for any $i$, and $\Gamma_i^+ \rightarrow \Gamma$, and
$\Gamma_i^- \rightarrow \Gamma$.

By Lemma 2.2, for any $\Gamma_i^+\subset \Si$, there exists a least area plane $\Sigma_i^+$ in $\BH^3$ asymptotic to
$\Gamma_i^+$. This defines a sequence of least area planes $\{\Sigma_i^+\}$. Now, by using Lemma 2.4, we take the limit
of a convergent subsequence. In the limit we get a collection of least area planes $\widehat{\Sigma}^+$ with
$\PI\widehat{\Sigma}^+ = \Gamma$, as $\PI\Sigma_i^+ = \Gamma_i^+ \rightarrow \Gamma$.

Now, we claim that the collection $\widehat{\Sigma}^+$ consists of only one least area plane. Assume that there are two
least area planes $\Sigma_a^+$ and $\Sigma_b^+$ in the collection $\widehat{\Sigma}^+$. Since $\PI \Sigma_a^+ = \PI
\Sigma_b^+ = \Gamma$, $\Sigma_a^+$ and $\Sigma_b^+$ might not be disjoint, but they are disjoint from least area planes
in the sequence, i.e. $\Sigma_i^+\cap \Sigma_{a,b}^+ =\emptyset$ for any $i$, by Lemma 3.1.

If $\Sigma_a^+$ and $\Sigma_b^+$ are disjoint, say $\Sigma_a^+$ is {\em above} $\Sigma_b^+$. By Lemma 3.1, we know that
for any $i$, $\Sigma_i^+$ is {\em above} both $\Sigma_a^+$ and $\Sigma_b^+$. However this means that $\Sigma_a^+$ is a
barrier between the sequence $\{\Sigma_i^+\}$ and $\Sigma_b^+$, and so, $\Sigma_b^+$ cannot be limit of this sequence,
which is a contradiction.

If $\Sigma_a^+$ and $\Sigma_b^+$ are not disjoint, then they intersect each other, and in some region, $\Sigma_b^+$ is
{\em above} $\Sigma_a^+$. However since $\Sigma_a^+$ is the limit of the sequence $\{\Sigma_i^+\}$, this would imply
$\Sigma_b^+$ must intersect planes $\Sigma_i^+$ for sufficiently large $i$. However, this contradicts the fact that
$\Sigma_b^+$ is disjoint from $\Sigma_i^+$ for any $i$, as they have disjoint asymptotic boundary. So, there exists a
unique least area plane $\Sigma^+$ in the collection $\widehat{\Sigma}^+$. Similarly, $\widehat{\Sigma}^- = \Sigma^-$.
By using similar arguments, one can conclude that these least area planes $\Sigma^+$, and $\Sigma^-$ are canonical,
i.e. independent of the choice of the sequence $\{\Gamma_i^\pm\}$ and $\{\Sigma_i^\pm\}$.

Now, let $\Sigma '$ be any least area plane with $\PI\Sigma '=\Gamma$. If $\Sigma '\cap\Sigma^+ \neq \emptyset$, then
some part of $\Sigma '$ must be {\em above} $\Sigma^+$. Since $\Sigma^+ =\lim \Sigma_i^+$, for sufficiently large $i$,
$\Sigma ' \cap\Sigma_i^+ \neq\emptyset$. However, $\PI\Sigma_i^+ = \Gamma_i^+$ is disjoint from $\Gamma=\PI\Sigma '$.
Then, by Lemma 3.1, $\Sigma '$ must be disjoint from $\Sigma_i^+$. This is a contradiction.

Similarly, this is true for $\Sigma^-$, too. Moreover, let $N \subset \BH^3$ be the region between $\Sigma^+$ and
$\Sigma^-$, i.e. $\partial N= \Sigma^+ \cup \Sigma^-$. Then by construction,  $N$ is also a canonical region for
$\Gamma$, and for any least area plane $\Sigma '$ with $\PI \Sigma '=\Gamma$, $\Sigma '$ is contained in the region
$N$, i.e. $\Sigma '\subset N$. This shows that if $\Sigma^+ = \Sigma^-$, there exists a unique least area plane
asymptotic to $\Gamma$. If $\Sigma^+ \neq \Sigma^-$, then they must be disjoint.
\end{pf}

Now, we are going to prove the main theorem of the section. This theorem says that for a generic simple closed curve in
$\Si$, there exists a unique least area plane.

\begin{thm} Let $A$ be the space of simple closed curves in $\Si$ and let $A'\subset A$ be the subspace containing the simple closed
curves in $\Si$ bounding a unique least area plane in $\BH^3$. Then, $A'$ is generic in $A$, i.e. $A-A'$ is a set of
first category.
\end{thm}

\begin{pf} We will prove this theorem in 2 steps.

\textbf{Claim 1:} $A'$ is dense in $A$ as a subspace of $C^0(S^1,\Si)$ with the supremum metric.\\

\begin{pf} $A$ is the space of Jordan curves in $\Si$. Then, $A= \{\alpha\in C^0(S^1,S^2)\ | \ \alpha(S^1) \mbox{ is an
embedding}\}$.

Now, let $\Gamma_0\in A$ be a simple closed curve in $\Si$. Since $\Gamma_0$ is simple, there exists a small
neighborhood $N(\Gamma_0)$ of $\Gamma_0$ which is an annulus in $\Si$. Let $\Gamma:(-\epsilon,\epsilon)\rightarrow A$
be a small path in $A$ through $\Gamma_0$ such that $\Gamma(t)=\Gamma_t$ and $\{\Gamma_t\}$ foliates $N(\Gamma)$ with
simple closed curves $\Gamma_t$. In other words, $\{\Gamma_t\}$ are pairwise disjoint simple closed curves, and
$N(\Gamma_0)=\bigcup_{t\in (-\epsilon,\epsilon)} \Gamma_t$.

Now, $N(\Gamma_0)$ separates $\Si$ into two parts, say $D^+$ and $D^-$, i.e. $\Si=N(\Gamma_0)\cup D^+\cup D^-$. Let
$p^+$ be a point in $D^+$ and let $p^-$ be a point in $D^-$ such that for a small $\delta$, $B_\delta(p^\pm)$ are in
the interior of $D^\pm$. Let $\beta$ be the geodesic in $\BH^3$ asymptotic to $p^+$ and $p^-$.

\begin{figure}[t]
{\epsfxsize=3.5in

\centerline{\epsfbox{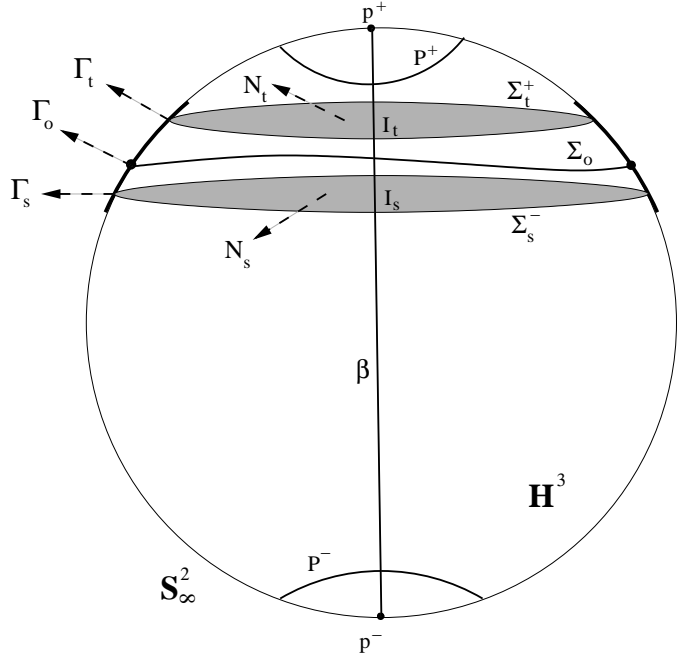}}}

\caption{\label{fig:figure1} {A finite segment of geodesic $\gamma$ intersects the collection of least area planes
$\Sigma_t$ in $\BH^n$ asymptotic to $\Gamma_t$ in $\SI$.}}
\end{figure}

By Lemma 3.2, for any $\Gamma_t$ either there exists a unique least area plane $\Sigma_t$ in $\BH^3$, or there is a
canonical region $N_t$ in $\BH^3$ between the canonical least area planes $\Sigma_t^+$ and $\Sigma_t^-$ (In Figure 1,
$\Gamma_t$ and $\Gamma_s$ bound more than one least area plane in $\BH^3$, whereas $\Gamma_0$ bounds a unique least
area plane $\Sigma_0$ in $\BH^3$). With abuse of notation, if $\Gamma_t$ bounds a unique least area plane $\Sigma_t$ in
$\BH^3$, define $N_t=\Sigma_t$ as a degenerate canonical neighborhood for $\Gamma_t$ (In Figure 1, $N_t$ and $N_s$
represent nondegenerate canonical neighborhoods, and $N_0=\Sigma_0$ represents degenerate canonical neighborhood.).
Then, let $\widehat{N}= \{N_t\}$ be the collection of these degenerate and nondegenerate canonical neighborhoods for
$t\in(-\epsilon,\epsilon)$. Clearly, degenerate neighborhood $N_t$ means $\Gamma_t$ bounds unique least area plane, and
nondegenerate neighborhood $N_s$ means that $\Gamma_s$ bounds more than one least area plane. Note that by Lemma 3.1,
all canonical neighborhoods in the collection are pairwise disjoint. On the other hand, by construction the geodesic
$\beta$ intersects all the canonical neighborhoods in the collection $\widehat{N}$.

We claim that the part of $\beta$ which intersects $\widehat{N}$ is a finite line segment. Let $P^+$ be the
geodesic plane asymptotic to round circle $\partial B_\delta(p^+)$ in $D^+$. Similarly, define $P^-$. By
Lemma 3.1, $P^\pm$ are disjoint from the collection of canonical regions $\widehat{N}$. Let $\beta\cap
P^\pm=\{q^\pm\}$. Then the part of $\beta$ which intersects $\widehat{N}$ is the line segment $l\subset
\beta$ with endpoints $q^+$ and $q^-$. Let $C$ be the length of this line segment $l$.

Now, for each $t\in(-\epsilon,\epsilon)$, we will assign a real number $s_t\geq 0$. If there exists a unique least area
plane $\Sigma_t$ in $\BH^3$ for $\Gamma_t$ ($N_t$ is degenerate), then let $s_t$ be $0$. If not, let $I_t = \beta\cap
N_t$, and $s_t$ be the length of $I_t$. Clearly if $\Gamma_t$ bounds more than one least area plane ($N_t$ is
nondegenerate), then $s_t > 0$. Also, it is clear that for any $t$, $I_t\subset l$ and $I_t\cap I_s=\emptyset$ for any
$t\neq s$. Then, $\sum_{t\in(-\epsilon,\epsilon)} s_t < C$ where $C$ is the length of $l$. This means for only
countably many $t\in(-\epsilon,\epsilon)$, $s_t > 0$. So, there are only countably many nondegenerate $N_t$ for
$t\in(-\epsilon,\epsilon)$. Hence, for all other $t$, $N_t$ is degenerate. This means there exist uncountably many
$t\in(-\epsilon,\epsilon)$, where $\Gamma_t$ bounds a unique least area plane. Since $\Gamma_0$ is arbitrary, this
proves $A '$ is dense in $A$.
\end{pf}

\textbf{Claim 2:} $A'$ is generic in $A$, i.e. $A-A'$ is a set of first category.\\

\begin{pf} We will prove that $A '$ is countable intersection of open dense subsets of a complete metric space. Then the
result will follow by Baire category theorem.

Since the space of continuous maps from circle to sphere $C^0(S^1,S^2)$ is complete with supremum metric, then
the closure of $A$ in $C^0(S^1,S^2)$, $\bar{A}\subset C^0(S^1,S^2)$, is also complete.

Now, we will define a sequence of open dense subsets $U^i\subset A$ such that their intersection will give us $A '$.
Let $\Gamma\in A$ be a simple closed curve in $\Si$, as in the Claim 1. Let $N(\Gamma)\subset \Si$ be a neighborhood of
$\Gamma$ in $\Si$, which is an open annulus. Then, define an open neighborhood $U_\Gamma$ of $\Gamma$ in $A$, such that
$U_\Gamma = \{\alpha \in A \ | \ \alpha(S^1)\subset N(\Gamma), \ \alpha \mbox{ is homotopic to } \Gamma\}$. Clearly,
$A= \bigcup_{\Gamma\in A} U_\Gamma$.  Now, define a geodesic $\beta_\Gamma$ as in Claim 1, which intersects all the
least area planes asymptotic to curves in $U_\Gamma$.

Now, for any $\alpha \in U_\Gamma$, by Lemma 3.2, there exists a canonical region $N_\alpha$ in $\BH^3$ (which can be
degenerate if $\alpha$ bounds a unique least area plane). Let $I_{\alpha,\Gamma} = N_\alpha \cap \beta_\Gamma$. Then
let $s_{\alpha,\Gamma}$ be the length of $I_{\alpha,\Gamma}$ ($s_{\alpha,\Gamma}$ is $0$ if $N_\alpha$ degenerate).
Hence, for every element $\alpha$ in $U_\Gamma$, we assign a real number $s_{\alpha,\Gamma} \geq 0$.

Now, we define the sequence of open dense subsets in $U_\Gamma$. Let $U^i_\Gamma = \{\alpha\in U_\Gamma \ | \
s_{\alpha,\Gamma} < 1/i \  \}$. We claim that $U^i_\Gamma$ is an open subset of $U_\Gamma$ and $A$. Let $\alpha\in
U^i_\Gamma$, and let $s_{\alpha,\Gamma} = \lambda < 1/i$. So, the interval $I_{\alpha,\Gamma}\subset \beta_\Gamma$ has
length $\lambda$. let $I ' \subset \beta_\Gamma$ be an interval containing $I_{\alpha,\Gamma}$ in its interior, and has
length less than $1/i$. By the proof of Claim 1, we can find two simple closed curves $\alpha^+, \alpha^- \in U_\Gamma$
with the following properties.

\begin{itemize}

\item $\alpha^\pm$ are disjoint from $\alpha$,

\item $\alpha^\pm$ are lying in opposite sides of $\alpha$ in $\Si$,

\item $\alpha^\pm$ bounds unique least area planes $\Sigma_{\alpha^\pm}$,

\item $\Sigma_{\alpha^\pm} \cap \beta_\Gamma \subset I '$.

\end{itemize}

The existence of such curves is clear from the proof Claim 1, as if one takes any foliation $\{\alpha_t\}$ of a small
neighborhood of $\alpha$ in $\Si$, there are uncountably many curves in the family bounding a unique least area plane,
and one can choose sufficiently close pair of curves to $\alpha$, to ensure the conditions above.

After finding $\alpha^\pm$, consider the open annulus $F_\alpha$ in $\Si$ bounded by $\alpha^+$ and
$\alpha^-$. Let $V_\alpha = \{ \gamma\in U_\Gamma \ | \ \gamma(S^1)\subset F_\alpha , \ \gamma \mbox{ is
homotopic to } \alpha \}$. Clearly, $V_\alpha$ is an open subset of $U_\Gamma$. If we can show
$V_\alpha\subset U^i_\Gamma$, then this proves $U^i_\Gamma$ is open for any $i$ and any $\Gamma\in A$.

Let $\gamma\in V_\alpha$ be any curve, and $N_\gamma$ be its canonical neighborhood given by Lemma 3.2. Since
$\gamma(S^1)\subset F_\alpha$, $\alpha^+$ and $\alpha^-$ lie in opposite sides of $\gamma$ in $\Si$. This
means $\Sigma_{\alpha^+}$  and $\Sigma_{\alpha^-}$ lie in opposite sides of $N_\gamma$. By choice of
$\alpha^\pm$, this implies $N_\gamma \cap \beta_\Gamma= I_{\gamma,\Gamma} \subset I '$. So, the length
$s_{\gamma,\Gamma}$ is less than $1/i$. This implies $\gamma\in U^i_\Gamma$, and so $V_\alpha\subset
U^i_\Gamma$. Hence, $U^i_\Gamma$ is open in $U_\Gamma$ and $A$.

Now, we can define the sequence of open dense subsets. let $U^i = \bigcup_{\Gamma\in A} U^i_\Gamma$ be an open subset
of $A$. Since, the elements in $A '$ represent the curves bounding a unique least area plane, for any $\alpha\in A '$,
and for any $\Gamma\in A$, $s_{\alpha,\Gamma} = 0$. This means $A'\subset U^i$ for any $i$. By Claim 1, $U^i$ is open
dense in $A$ for any $i>0$.

As we mention at the beginning of the proof, since the space of continuous maps from circle to sphere $C^0(S^1,S^2)$ is
complete with supremum metric, then the closure $\bar{A}$ of $A$ in $C^0(S^1,S^2)$ is also complete metric space. Since
$A'$ is dense in $A$, it is also dense in $\bar{A}$. As $A$ is open in $C^0(S^1,S^2)$, this implies $U^i$ is a sequence
of open dense subsets of $\bar{A}$. On the other hand, since $s_{\alpha,\Gamma} = 0$ for any $\alpha\in A '$, and for
any $\Gamma\in A$, $A ' = \bigcap_{i>0} U^i$. Then, $A-A'$ is a set of first category, by Baire Category Theorem.
Hence, $A'$ is generic in $A$.
\end{pf}
\end{pf}

\begin{rmk} This result is similar to the generic uniqueness result in \cite{Co2}. In \cite{Co2}, we used a heavy
machinery of analysis to prove that there exists an open dense subset in the space of $C^{3,\mu}$-smooth embeddings of
circle into sphere, where any simple closed curve in this space bounds a unique least area plane in $\BH^3$. In the
above result, the argument is fairly simple, and does not use the analytical machinery.
\end{rmk}

\section{Area Minimizing Hypersurfaces in $\BH^n$}

In this section, we will show that the space of codimension-$1$ closed submanifolds of $\SI$ bounding a unique
absolutely area minimizing hypersurface in $\BH^n$ is dense in the space of all codimension-$1$ closed submanifolds of
$\SI$. Indeed, we will show that in some sense, a generic closed manifold $\Gamma$ in $\SI$ bounds a unique absolutely
area minimizing hypersurface $\Sigma$ in $\BH^n$. The idea is similar to the previous section.

First, we need to show a simple topological lemma.

\begin{lem}
Any codimension-$1$ closed submanifold $\Gamma$ of $S^{n-1}$ is separating. Moreover, if $\Sigma$ is a hypersurface
with boundary in the closed unit ball $B^n$, where the boundary $\partial \Sigma \subset
\partial B^n$, then $\Sigma$ is also separating, too.
\end{lem}

\begin{pf}
A non-separating codimension-$1$ closed submanifold in $S^{n-1}$ gives a nontrivial homology in $n-2$ level. However
since $H_{n-2}(S^{n-1})$ is trivial, this is a contradiction.

Let $\Sigma$ be as in the assumption. Take the double of $B^n$, then $B^n \sqcup \widehat{B^n} =S^n$, and
$\Sigma \sqcup \widehat{\Sigma}$ is a codimension-$1$ closed submanifold of $S^n$. By above, $\Sigma \sqcup
\widehat{\Sigma}$ is separating in $S^n$. Hence, $\Sigma$ is separating in $B^n$.
\end{pf}

Now, we will prove a disjointness lemma analogous to Lemma 3.1. This lemma roughly says that if asymptotic boundaries
of two absolutely area minimizing hypersurfaces in $\BH^n$ are disjoint in $\SI$, then they are disjoint in $\BH^n$.

\begin{lem}
Let $\Gamma_1$ and $\Gamma_2$ be two disjoint connected closed codimension-$1$ submanifolds in $\SI$. If
$\Sigma_1$ and $\Sigma_2$ are absolutely area minimizing hypersurfaces in $\BH^n$ with $\PI \Sigma_i =
\Gamma_i$, then $\Sigma_1$ and $\Sigma_2$ are disjoint, too.
\end{lem}

\begin{pf}
Assume that the absolutely area minimizing hypersurfaces are not disjoint, i.e. $\Sigma_1\cap\Sigma_2\neq\emptyset$. By
Lemma 4.1, $\Sigma_1$, and $\Sigma_2$ separates $\BH^n$ into two parts. So, say $\BH^n-\Sigma_i =
\Omega^+_i\cup\Omega^-_i$.

Now, consider the intersection of hypersurfaces $\alpha=\Sigma_1\cap\Sigma_2$. Since the asymptotic boundaries
$\Gamma_1$ and $\Gamma_2$ are disjoint in $\SI$, by using the regularity results of Hardt and Lin in \cite{HL}, we can
conclude that the intersection set $\alpha$ is in compact part of $\BH^n$. Moreover, by maximum principle \cite{Si},
the intersection cannot have isolated tangential intersections.

\begin{figure}[t]
{\epsfxsize=5in

\centerline{\epsfbox{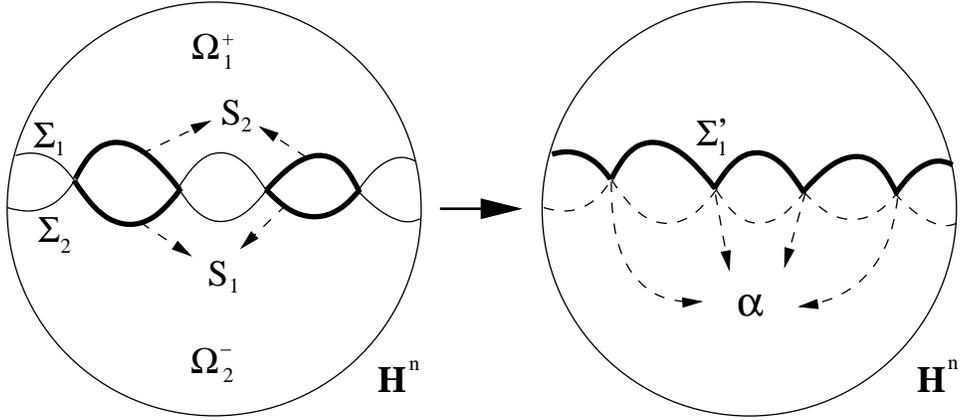}}}

\caption{\label{fig:figure1} {$S_1$ is the part of $\Sigma_1$ lying below $\Sigma_2$, and $S_2$ is the part of
$\Sigma_2$ lying above $\Sigma_1$. After swaping $S_1$ and $S_2$, we get a new area minimizing hypersurface $\Sigma_1
'$ with singularity along $\alpha=\Sigma_1\cap\Sigma_2$.}}
\end{figure}

Now, without loss of generality, we assume that $\Sigma_1$ is {\em above} $\Sigma_2$ (the noncompact part of $\Sigma_1$
lies in $\Omega^+_2$). Now define the compact subhypersurfaces $S_i$ in $\Sigma_i$ as $S_1=\Sigma_1\cap \Omega^-_2$,
and $S_2 = \Sigma_2\cap\Omega^+_1$.  In other words, $S_1$ is the part of $\Sigma_1$ lying {\em below} $\Sigma_2$, and
$S_2$ is the part of $\Sigma_2$ lying {\em above} $\Sigma_1$. Then, $\partial S_1 =\partial S_2 =\alpha$. See Figure 2.

On the other hand, since $\Sigma_1$ and $\Sigma_2$ are absolutely area minimizing, then by definition, so are $S_1$ and
$S_2$, too. Then by swaping the surfaces, we can get new absolutely area minimizing hypersurfaces. In other words, let
$\Sigma_1 ' = \{\Sigma_1-S_1\}\cup S_2$, and $\Sigma_2 ' = \{\Sigma_2-S_2\}\cup S_1$ are also absolutely area
minimizing hypersurfaces. However, in this new hypersurfaces, we will have a singularity set along $\alpha$, which
contradicts to regularity theorem for absolutely area minimizing hypersurfaces, i.e. Lemma 2.3.
\end{pf}

\begin{lem}
Let $\Gamma$ be a connected closed codimension-$1$ submanifold of $\SI$. Then either there exists a unique
absolutely area minimizing hypersurface $\Sigma$ in $\BH^n$ asymptotic to $\Gamma$, or there are two canonical
disjoint extremal absolutely area minimizing hypersurfaces $\Sigma^+$ and $\Sigma^-$ in $\BH^n$ asymptotic to
$\Gamma$.
\end{lem}

\begin{pf}
Let $\Gamma$ be a connected closed codimension-$1$ submanifold of $\SI$. Then by Lemma 4.1, $\Gamma$ separates $\SI$
into two parts, say $\Omega^+$ and $\Omega^-$. Define sequences of pairwise disjoint closed submanifolds of the same
topological type $\{\Gamma_i^+\}$ and $\{\Gamma_i^-\}$ in $\SI$ such that $\Gamma_i^+\subset \Omega^+$, and
$\Gamma_i^-\subset \Omega^-$ for any $i$, and $\Gamma_i^+ \rightarrow \Gamma$, and $\Gamma_i^- \rightarrow \Gamma$ in
Hausdorff metric. In other words, $\{\Gamma_i^+\}$ and $\{\Gamma_i^-\}$ converges to $\Gamma$ from opposite sides.

By Lemma 2.1, for any $\Gamma_i^+\subset \Si$, there exists an absolutely area minimizing hypersurface $\Sigma_i^+$ in
$\BH^n$. This defines a sequence of absolutely area minimizing hypersurfaces $\{\Sigma_i^+\}$. By Lemma 2.5, we get a
convergent subsequence $\Sigma_{i_j}^+\rightarrow \Sigma^+$. Hence, we get the absolutely area minimizing hypersurface
$\Sigma^+$ in $\BH^n$ asymptotic to $\Gamma$. Similarly, we get the absolutely area minimizing hypersurface $\Sigma^-$
in $\BH^n$ asymptotic to $\Gamma$. Similar arguments show that these absolutely area minimizing hypersurfaces
$\Sigma^\pm$ are canonical by their construction, i.e. independent of the choice of the sequence $\{\Gamma_i^\pm\}$ and
$\{\Sigma_i^\pm\}$.

Assume that $\Sigma^+ \neq\Sigma^-$, and they are not disjoint. Since these are absolutely area minimizing
hypersurfaces, nontrivial intersection implies some part of $\Sigma^-$ lies {\em above} $\Sigma^+$. Since
$\Sigma^+=\lim \Sigma_{i_j}^+$, $\Sigma^-$ must also intersect some $\Sigma_{i_j}^+$ for sufficiently large $i_j$.
However by Lemma 4.2, $\Sigma_{i_j}^+$ is disjoint from $\Sigma^-$ as $\PI\Sigma_{i_j}^+ = \Gamma_{i_j}^+$ is disjoint
from $\PI \Sigma^- = \Gamma$. This is a contradiction. This shows $\Sigma^+$ and $\Sigma^-$ are disjoint.

Similar arguments show that $\Sigma^\pm$ are disjoint from any absolutely area minimizing hypersurface $\Sigma '$
asymptotic $\Gamma$. As the sequences of $\Sigma_i^+$ and $\Sigma_i^-$ forms a barrier for other absolutely area
minimizing hypersurfaces asymptotic to $\Gamma$, any such absolutely area minimizing hypersurface must lie in the
region bounded by $\Sigma^+$ and $\Sigma^-$ in $\BH^n$. This shows that if $\Sigma^+ = \Sigma^-$, then there exists a
unique absolutely area minimizing hypersurface asymptotic to $\Gamma$.
\end{pf}

\begin{rmk}
By above theorem and its proof, if $\Gamma$ bounds more than one absolutely area minimizing hypersurface, then there
exists a canonical region $N_\Gamma$ in $\BH^n$ asymptotic to $\Gamma$ such that $N_\Gamma$ is the region between the
canonical absolutely area minimizing hypersurfaces $\Sigma^+$ and $\Sigma^-$. Moreover, by using similar ideas to the
proof of Lemma 3.2, one can show that any absolutely area minimizing hypersurface in $\BH^n$ asymptotic to $\Gamma$ is
in the region $N_\Gamma$.
\end{rmk}

Now, we can prove the main result of the paper.

\begin{thm}
Let $B$ be the space of connected closed codimension-$1$ submanifolds of $\SI$, and let $B'\subset B$ be the
subspace containing the closed submanifolds of $\SI$ bounding a unique absolutely area minimizing hypersurface
in $\BH^n$. Then $B'$ is dense in $B$.
\end{thm}

\begin{pf}
Let $B$ be the space of connected closed codimension-$1$ submanifolds of $\SI$ with Hausdorff metric. Let $\Gamma_0\in
B$ be a closed submanifold in $\SI$. Since $\Gamma_0$ is closed submanifold, there exists a small regular neighborhood
$N(\Gamma_0)$ of $\Gamma_0$ in $\SI$, which is homeomorphic to $\Gamma_0\times I$. Let
$\Gamma:(-\epsilon,\epsilon)\rightarrow B$ be a small path in $B$ through $\Gamma_0$ such that $\Gamma(t)=\Gamma_t$ and
$\{\Gamma_t\}$ foliates $N(\Gamma_0)$ with closed submanifolds homeomorphic to $\Gamma_0$. In other words,
$\{\Gamma_t\}$ are pairwise disjoint closed submanifolds homeomorphic to $\Gamma_0$, and $N(\Gamma_0)=\bigcup_{t\in
(-\epsilon,\epsilon)} \Gamma_t$.

By Lemma 4.1, $N(\Gamma_0)$ separates $\SI$ into two parts, say $\Omega^+$ and $\Omega^-$, i.e. $\SI=N(\Gamma_0)\cup
\Omega^+\cup \Omega^-$. Let $p^+$ be a point in $\Omega^+$ and let $p^-$ be a point in $\Omega^-$ such that for a small
$\delta$, $B_\delta(p^\pm)$ are in the interior of $\Omega^\pm$. Let $\beta$ be the geodesic in $\BH^n$ asymptotic to
$p^+$ and $p^-$.

By Lemma 4.3 and Remark 4.1, for any $\Gamma_t$ either there exists a unique absolutely area minimizing hypersurface
$\Sigma_t$ in $\BH^n$, or there is a canonical region $N_t$ in $\BH^n$ asymptotic to $\Gamma_t$, namely the region
between the canonical absolutely area minimizing hypersurfaces $\Sigma_t^+$ and $\Sigma_t^-$. With abuse of notation,
if $\Gamma_t$ bounds a unique absolutely area minimizing hypersurface $\Sigma_t$ in $\BH^n$, define $N_t=\Sigma_t$ as a
degenerate canonical neighborhood for $\Gamma_t$. Then, let $\widehat{N}= \{N_t\}$ be the collection of these
degenerate and nondegenerate canonical neighborhoods for $t\in(-\epsilon,\epsilon)$. Clearly, degenerate neighborhood
$N_t$ means $\Gamma_t$ bounds a unique absolutely area minimizing hypersurface, and nondegenerate neighborhood $N_s$
means that $\Gamma_s$ bounds more than one absolutely area minimizing hypersurfaces. Note that by Lemma 4.2, all
canonical neighborhoods in the collection are pairwise disjoint. On the other hand, by Lemma 4.1, the geodesic $\beta$
intersects all the canonical neighborhoods in the collection $\widehat{N}$.

We claim that the part of $\beta$ which intersects $\widehat{N}$ is a finite line segment. Let $P^+$ be the
geodesic hyperplane asymptotic to round sphere $\partial B_\delta(p^+)$ in $\Omega^+$. Similarly, define
$P^-$. By Lemma 4.2, $P^\pm$ are disjoint from the collection of canonical regions $\widehat{N}$. Let
$\beta\cap P^\pm=\{q^\pm\}$. Then the part of $\beta$ which intersects $\widehat{N}$ is the line segment
$l\subset \beta$ with endpoints $q^+$ and $q^-$. Let $C$ be the length of this line segment $l$.

Now, for each $t\in(-\epsilon,\epsilon)$, we will assign a real number $s_t\geq 0$. If there exists a unique absolutely
area minimizing hypersurface $\Sigma_t$ for $\Gamma_t$ ($N_t$ is degenerate), then let $s_t$ be $0$. If not, let $I_t =
\beta\cap N_t$, and $s_t$ be the length of $I_t$. Clearly if $\Gamma_t$ bounds more than one least area plane ($N_t$ is
nondegenerate), then $s_t > 0$. Also, it is clear that for any $t$, $I_t\subset l$ and $I_t\cap I_s=\emptyset$ for any
$t\neq s$. Then, $\sum_{t\in(-\epsilon,\epsilon)} s_t < C$ where $C$ is the length of $l$. This means for only
countably many $t\in(-\epsilon,\epsilon)$, $s_t > 0$. So, there are only countably many nondegenerate $N_t$ for
$t\in(-\epsilon,\epsilon)$. Hence, for all other $t$, $N_t$ is degenerate. This means there exist uncountably many
$t\in(-\epsilon,\epsilon)$, where $\Gamma_t$ bounds a unique absolutely area minimizing hypersurface. Since $\Gamma_0$
is arbitrary, this proves $B'$ is dense in $B$.
\end{pf}

\begin{rmk}
This density result can be generalized to a genericity in some sense. Here, if one does not consider the whole space of
closed codimension-1 submanifolds of $\SI$, but specify the topological type of the closed submanifolds, by using
similar arguments to the proof of Claim 2 in Theorem 3.3, one can get a genericity result in that space. In other
words, one can stratify the whole space of codimension-1 closed submanifolds by topological type, and get the
genericity result in each strata.
\end{rmk}

\begin{cor}
Let $A$ be the space of simple closed curves in $\Si$ and let $A'\subset A$ be the subspace containing the simple
closed curves in $\Si$ bounding a unique absolutely area minimizing surface in $\BH^3$. Then, $A'$ is generic in $A$,
i.e. $A-A'$ is a set of first category.
\end{cor}

\begin{pf}
By Theorem 4.4, we know that $A'$ is dense in $A$. Then by using the proof of Claim 2 in Theorem 3.3 in this setting,
it is clear that $A'$ is generic in $A$ as a subspace of $C^0(S^1,S^2)$ with supremum metric. In other words, $A-A'$ is
a set of first category.
\end{pf}

\begin{rmk}
This corollary shows that like in the case of least area planes (Theorem 3.3), a generic simple closed curve in $\Si$
also bounds a unique absolutely area minimizing surface.
\end{rmk}

\section{Nonuniqueness results}

In this section, we will show that there exists a simple closed curve in $\Si$ which is the asymptotic boundary of more
than one absolutely area minimizing surface in $\BHH$. First, we need a lemma about the limits of absolutely area
minimizing surfaces.

A short outline of the method is the following. We first construct a simple closed curve $\Gamma_0$ in $\Si$ which is
the asymptotic boundary of more than one \textit{minimal surface} in $\BHH$. Then, we foliate $\Si$ with simple closed
curves $\{\Gamma_t\}$ where $\Gamma_0$ is a leaf in the foliation. We show that if each $\Gamma_t$ bounds a unique
absolutely area minimizing surface $\Sigma_t$ in $\BHH$, then the family of surfaces $\{\Sigma_t\}$ must foliate the
whole $\BHH$. However, since we chose $\Gamma_0$ to bound more than one minimal surface in $\BHH$, one of the surfaces
must have a tangential intersection with one of the leaves in the foliation. This contradicts to the maximum principle
for minimal surfaces.

Now, we quote a result on the existence of simple closed curves in $\Si$ which are the asymptotic boundaries of more
than one minimal surface in $\BH^3$.

\begin{figure}[b]

\relabelbox \small {\epsfxsize=5in
  \centerline{\epsfbox{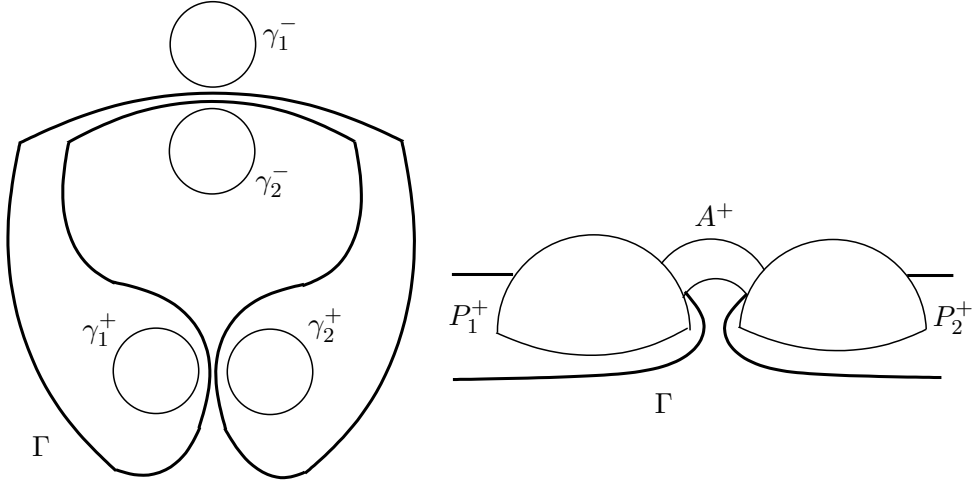}}}

\relabel{1}{{$\gamma_1^-$}}

\relabel{2}{$\gamma_2^-$}

\relabel{3}{$\gamma_1^+$}

\relabel{4}{$\gamma_2^+$}

\relabel{5}{$P_1^+$}

\relabel{6}{$P_2^+$}

\relabel{7}{$A^+$}

\relabel{8}{$\Gamma$}

\relabel{9}{$\Gamma$}

\endrelabelbox

\caption{\label{fig:figure1} {$\Gamma$ is a simple closed curve in $\Si$. $\gamma_i^+$ and $\gamma_{i}^-$ are round
circles in $\Si$ bounding the geodesic planes $P_i^+$ and $P_{i}^-$ in $\BHH$.}}

\end{figure}

\begin{lem} \cite{A2} There is a set $\Delta$ of Jordan curves in $\Si$ such that for any $\Gamma\in \Delta$, there
exist infinitely many complete, smoothly embedded minimal surfaces in $\BH^3$ asymptotic to $\Gamma$.
\end{lem}

\begin{rmk} Alternatively, one can construct simple closed curves in $\Si$ bounding more than one minimal surface in $\BH^3$ as
follows. By using the technique in [Ha], one can construct a simple closed curve $\Gamma$ in $\Si$ such that the
absolutely area minimizing surface $\Sigma$ asymptotic to $\Gamma$ has positive genus (See Figure 3). Then, $\Sigma$
separates $\BH^3$ into two parts $\Omega^+$ and $\Omega^-$ which are both mean convex domains. Then by using
Meeks-Yau's results in \cite{MY}, one can get sequences of least area disks $\{D_i^\pm\}$ in $\Omega^\pm$ whose
boundaries converges to $\Gamma$ in $\Si$. By taking the limit, one can get two least area planes to $\Gamma$. We used
$\Sigma$ as barrier to get distinct limit from the sequences $\{D_i^+\}$ and $\{D_i^-\}$. Note that the planes $M^\pm$
are least area just in $\Omega^\pm$. So, $M^\pm$ may not be least area in $\BH^3$, but of course, they are still
minimal planes.
\end{rmk}

\begin{thm} There exists a simple closed curve $\Gamma$ in $\Si$ such that $\Gamma$ bounds more than
one absolutely area minimizing surface $\{\Sigma_i\}$ in $\BHH$, i.e. $\PI \Sigma_i = \Gamma$.
\end{thm}

\begin{pf} Assume that for any simple closed curve $\Gamma$ in $\Si$, there exists a unique
complete absolutely area minimizing surface $\Sigma$ in $\BHH$ with $\PI \Sigma =\Gamma$. Let $\Gamma_0$ be a simple
closed curve in $\Si$ such that $M_1$ and $M_2$ are two distinct minimal surfaces in $\BHH$ asymptotic to $\Gamma$ as
in Lemma 5.1.

Now, foliate $\Si$ by simple closed curves $\{\Gamma_t\}$ where $\Gamma_0$ is a leave in the foliation. Note that,
there are only two singular leaves in the foliation which are points, and all other leaves are embedded simple closed
curves in $\Si$. By assumption, for any $\Gamma_t$, there exists a unique absolutely area minimizing surface $\Sigma_t$
in $\BHH$.

We claim that $\{\Sigma_t\}$ will be a foliation of $\BHH$. By Lemma 4.2, $\Sigma_t\cap\Sigma_s = \emptyset$ for any
$s\neq t$. Hence, the only way to fail to be a foliation for $\{\Sigma_t\}$ is to have a gap between two leaves.

Now, assume that there is a gap between the leaves $\{\Sigma_t \ | \ t>s\ \}$ and $\{\Sigma_t \ | \ t<s\ \}$. This
implies if we have sequences $\{\Gamma_{t^+_i}\}$ and $\{\Gamma_{t^-_i}\}$, where $t^+_i\rightarrow s$ from positive
side, and $t^-_i\rightarrow s$ from negative side, then one of the sequences $\{\Sigma_{t^+_i}\}$ and
$\{\Sigma_{t^-_i}\}$ has no subsequences converging to $\Sigma_s$, because of the existence of the gap between the
leaves (Recall that we assume there is a unique absolutely area minimizing surface $\Sigma_t$ with $\PI \Sigma_t =
\Gamma_t$). However, this contradicts to Lemma 2.5. So, $\{\Sigma_t\}$ will be a foliation of $\BHH$.

Now, by assumption, there are two distinct minimal surface $M_1$ and $M_2$ which are asymptotic to $\Gamma_0$. This
implies at least one of these surfaces is not a leaf of the foliation, and say $M_2$, must intersect the leaves in the
foliation nontrivially. Since $\{\Sigma_t\}$ foliates whole $\BHH$, this means $M_2$ must intersect tangentially (lying
in one side) one of the leaves in the foliation. However, this contradicts to the maximum principle. So, one of the
simple closed curves in $\{\Gamma_t\}$ must bound more than one absolutely area minimizing surface in $\BHH$. The proof
follows.
\end{pf}

By using similar ideas, one can prove an analogous theorem for least area planes in $\BHH$.

\begin{thm} There exists a simple closed curve $\Gamma$ in $\Si$ such that $\Gamma$ bounds more than
one least area plane $\{P_i\}$ in $\BHH$, i.e. $\PI P_i = \Gamma$.
\end{thm}

\begin{pf} The proof is completely analogous to the proof of previous theorem. Again, we start with
the same foliation of $\Si$ with simple closed curves $\{\Gamma_t\}$ containing $\Gamma_0$ which bounds more than one
complete minimal surface in $\BHH$.

Assume that $P_t$ is the unique least area plane in $\BHH$ with $\PI P_t = \Gamma_t$. Like in the proof of the theorem
above, by using Lemma 3.1, and Lemma 4.2, it can be showed that the family of least area planes $\{P_t\}$ foliates
$\BHH$. However, like in above theorem, $\Gamma_0$ bounds more than one area minimizing surface, and at least one of
them is not a leaf of the foliation, say $M_2$. Hence, $M_2$ must intersect tangentially (lying in one side) one of the
leaves in the foliation. Again, this contradicts to the maximum principle for minimal surfaces. The proof follows.
\end{pf}

\begin{rmk} The same proof may not work for area minimizing surfaces in a specified topological class.
The problem is that Lemma 2.5 may not be true for this case as the limiting surface might not be in the same
topological class.
\end{rmk}

\begin{rmk}As the introduction suggests, there is no known example of a simple closed curve of $\Si$ with
nonunique solution to the asymptotic Plateau problem. Unfortunately, the results above show the existence of such an
example, but they do not give one. The main problem to find such an example is the noncompactness of the objects. In
compact case, the quantitative data enables you to find such examples like baseball curve on a sphere, but in the
asymptotic case the techniques do not work because of the lack of the quantitative data. To overcome this problem, it
might be possible to employ the renormalized area defined by Alexakis-Mazzeo in \cite{AM} in order to construct an
explicit simple closed curve in $\Si$ bounding more than one absolutely area minimizing surface.
\end{rmk}

\section{Concluding Remarks}

In this paper, we showed that the space of closed, codimension-$1$ submanifolds of $\SI$ has a dense (and generic in a
sense) subspace of closed, codimension-$1$ submanifolds of $\SI$ bounding a unique absolutely area minimizing
hypersurface in $\BH^n$. As we discussed in the introduction, Anderson showed this result for closed submanifolds
bounding convex domains in $\SI$ in \cite{A1}. Then, Hardt and Lin generalized this result to closed submanifolds
bounding star shaped domains in $\SI$ in \cite{HL}. These were the only cases known so far. Hence, our result shows how
abundant they are. In dimension $3$, they are generic, and in the higher dimensions, they are dense.

The technique which we employ here is very general, and it applies to many different settings of Plateau problem. In
particular, it can naturally be generalized to the Gromov-Hadamard spaces which is studied by Lang in \cite{L}, and it
can be generalized to the mean convex domains with spherical boundary which is studied by Lin in \cite{Li}. In other
words, codimension-1 closed submanifolds in the boundary of these spaces generically bounds unique absolutely area
minimizing hypersurfaces. Generalizing this technique in the context of constant mean curvature hypersurfaces in
hyperbolic space also gives similar results. On the other hand, they can also be applied in Gromov hyperbolic 3-spaces
with cocompact metric where the author solved the asymptotic Plateau problem \cite{Co3}.

On the other hand, it was not known whether all closed codimension-$1$ submanifolds in $\SI$ have a unique solution to
the asymptotic Plateau problem or not. The only known results about nonuniqueness are also come from Anderson in
\cite{A2}. He constructs examples of simple closed curves in $\Si$ bounding more than one complete \textit{minimal}
surface in $\BH^3$. These examples are also area minimizing in their topological class. However, none of them are
absolutely area minimizing, i.e. a solution to the asymptotic Plateau problem. In Section 5, we prove the existence of
simple closed curves in $\Si$ with nonunique solution to asymptotic Plateau problem, and hence, give an answer for
dimension $3$. However, there is no result in higher dimensions yet. In other words, it is not known whether there
exist closed codimension-$1$ submanifolds in $\SI$ bounding more than one absolutely area minimizing hypersurfaces or
not for $n>3$.

It might be possible to extend the techniques in this paper to address the nonuniqueness question in higher dimensions.
It is possible to use the technique in Remark 5.1 (Figure 3), to get a codimension-$1$ sphere $\Gamma_0$ in $\SI$
bounding an absolutely area minimizing hypersurface $\Sigma_0$ in $\BH^n$ which is not a hyperplane. Like in the proof
of Theorem 5.2, by foliating $\SI$ by closed, codimension-$1$ submanifolds $\{\Gamma_t\}$, and by assuming uniqueness
of absolutely area minimizing hypersurfaces, one can get a foliation of $\BH^n$ by the absolutely area minimizing
hypersurfaces $\{\Sigma_t\}$ with $\PI\Sigma_t=\Gamma_t$. Let $\gamma$ be a $n-2$-sphere in $\Sigma_0$ which is "close"
to $\Gamma_0$. Then by using \cite{Wh2}, one can get a compact area minimizing hyperplane $M$ whose boundary is
$\gamma$. Hence by construction, $M$ cannot be in a leaf in the foliation $\{\Sigma_t\}$. Like in the proof of Theorem
5.2, it might be possible to get a contradiction by studying the intersection of $M$ with the foliation $\{\Sigma_t\}$.

Unfortunately, the maximum principle \cite{Si} does not work here since while $\Sigma_t$ is absolutely area minimizing,
$M$ is not. Also, the maximum principles due to Solomon-White \cite{SW} and Ilmanen \cite{I} which are for stationary
varifolds are not enough to get a contradiction since the minimizing hyperplane $M$ given by \cite{Wh2} might have
codimension-$1$ singularities, while \cite{SW} and \cite{I} works up to codimension-$2$ singularities. So, to get a
contradiction here, one needs a stronger maximum principle, or more regular area minimizing hyperplane $M$.

\end{document}